# PROBABILITY AND STATISTICS IN GREECE: PAST, PRESENT, AND FUTURE*

Andreas N. Philippou[1], Apostolos Kasapis[2], and Alex Karagrigoriou[3]

[1]University of Patras, [2]Ministry of Economy and Finance, [3]University of Piraeus

## ABSTRACT

The history of official statistics in modern Greece dates to 1825, when Law No. 48 provided for census-based conscription. The first Population Census was conducted in 1828. A dedicated Department of Statistics was established within the Bureau of Public Economy in 1861, marking the first specialized statistical service in the Greek public administration. In 1925, the General Statistical Service of Greece was established to consolidate state statistical activities. An initial framework for its reorganization was introduced in 1953, and the National Statistical Service of Greece was established in 1956. It became a General Secretariat in 1986 and was transformed into the independent Hellenic Statistical Authority (ELSTAT) in 2010. Probability and statistics emerged as independent academic disciplines in Greece between the late 1960s and 1982. Between 1968 and 1980, seven chairs in probability, statistics or closely related fields were established at five Greek higher education institutions, including the Higher School of Economics and Commercial Sciences, which was renamed the Athens University of Economics and Business in 1989. Six of these chairs were filled. Following the abolition of the chair system in 1982, more than 200 academics with doctorates in probability, statistics, or related fields have served in Greek public universities. Probability, statistics, and related disciplines are now well established in Greece, and Greek contributions to these fields have gained international recognition.

## 1. Traces of Statistics and the First Census

In 2025, Greece marked two centuries since the enactment of the first law providing for population enumeration and one century since the establishment of the General Statistical Service of Greece. Traces of statistical thinking in Greece can be found in ancient medical and political writings. The Hippocratic Corpus contains some of the earliest surviving attempts to systematically observe and document the course of individual illnesses, while

* Based on an invited talk delivered at the ΣtatMOD2025 Conference, University of Piraeus, Greece. Unless otherwise indicated, publication counts were obtained from MathSciNet (https://mathscinet.ams.org/mathscinet) and cross-checked against zbMATH Open (https://zbmath.org/). Information on PhD students and academic descendants was drawn from the Mathematics Genealogy Project (https://www.mathgenealogy.org/) and, where necessary, supplemented by published biographical sources and information provided by the authors. All three databases were last accessed in July 2026.

Aristotle's *Politics* offers one of antiquity's most extensive examinations of population size, fertility and the composition of political communities. These writings reveal an early concern with systematic observation, classification and population-related reasoning.

Following the Greek War of Independence and the establishment of the modern Greek state, the first Population Census was conducted in 1828 under Ioannis Kapodistrias, the first Governor of modern Greece. This work built upon Law No. 48 of 10 September 1825, "On Census-Based Conscription" (Academy of Athens, Research Centre for the History of Greek Law, 1825). With Panoutsos Notaras presiding over the Legislative Body and Georgios Kountouriotis heading the Executive Body, the law provided for a ratio of one soldier per 100 inhabitants in each city, town or village. The 1828 Census recorded 753,400 inhabitants and also provided a retrospective estimate of 938,765 inhabitants for 1821 (Hellenic Statistical Authority, 2021).

## 2. The Early Institutional Development of Official Statistics (1833–1861)

Following the establishment of the Secretariat of the Interior in 1833, the Bureau of Public Economy was created within it by Royal Decree of 29 April 1834. Among its broad responsibilities was the development of statistical work. The Bureau subsequently conducted population censuses annually between 1836 and 1845 and again at wider intervals in 1848, 1853 and 1856, indicating that population enumeration had already become a recurring state activity before the establishment of a dedicated Department of Statistics in 1861 (Hellenic Statistical Authority, 2021). At the same time, other public services collected statistics within their respective fields. For example, the Ministry of Finance began compiling foreign trade statistics in 1853; these were initially published in the Government Gazette and, from 1859 onwards, in a separate publication on foreign trade and shipping (Hellenic Statistical Authority, 2021).

The Bureau of Public Economy was reorganized by Law FKD of 1859 and, under a Decree issued in 1861, was divided into four departments, one of which was the Department of Statistics. This marked the first establishment of a specialized statistical service within the Greek public administration. Its responsibilities included the population census, statistics on arts and trades, agriculture, shipping and commerce, and the broader collection of statistical information (Hellenic Statistical Authority, 2021).

## 3. The General Statistical Service of Greece (1925–1955), the National Statistical Service of Greece (1956–2009), and the Hellenic Statistical Authority (2010–Present)

### 3.1. Official Statistics

As discussed in Sections 1 and 2, the institutional development of official statistics in Greece began during the early years of the modern Greek state. The first Population Census was conducted in 1828, while a dedicated Department of Statistics was established in 1861.

By Law 3724/1910, a Department of Statistics was established within the Ministry of Agriculture, Commerce and Industry. The Ministry was renamed the Ministry of National Economy in 1911, and a Royal Decree subsequently defined the organization and operation of its statistical service. In 1913, the Department was upgraded to a Directorate of Statistics, responsible for the collection, verification, processing and publication of national statistics (Hellenic Statistical Authority, 2021).

By the Legislative Decree of 31 March 1925, the General Statistical Service of Greece (GSSG) was established within the Ministry of National Economy, consolidating statistical activities that had previously been dispersed across different ministries. The reform sought to ensure methodological uniformity, promote staff specialization, prevent the duplication of statistical work and improve the efficiency of statistical production. The GSSG subsequently conducted major population, agricultural and industrial censuses and issued regular publications, including the Monthly Statistical Bulletin and the Statistical Yearbook of Greece (Hellenic Statistical Authority, 2021).

Legislative Decree 2516/1953 provided an initial framework for the reorganization of the statistical service. The National Statistical Service of Greece (NSSG) was subsequently established under the Ministry of Coordination by Legislative Decree 3627/1956, which substantially amended and replaced the 1953 decree. Its organization was completed by Royal Decree 871/1960. In 1986, the General Secretariat of the National Statistical Service of Greece was established under the Ministry of National Economy (Hellenic Statistical Authority, 2021). Petros Kouvelis headed the NSSG from 1956 to 1970, while Emmanouil Kontopyrakis served as the last General Secretary of the NSSG from 2004 to 2009.

In 2010, Law 3832/2010 transformed the General Secretariat of the National Statistical Service of Greece into an independent administrative authority under the name Hellenic Statistical Authority (ELSTAT). The same law established the Hellenic Statistical System, comprising the rules, activities and agencies responsible for the development, production and dissemination of the country's official statistics (Hellenic Statistical Authority, 2021). The first President of ELSTAT was Andreas Georgiou, who served from 2010 to 2015. The Authority is currently headed by Athanassios Katsis, who assumed office in May 2026 (Hellenic Statistical Authority, 2026).

### 3.2. Statistical Production and Methodological Development

The methodological evolution of population censuses in Greece reflects a gradual transition from administrative enumeration to scientific statistical inquiry. Although the 1828 Population Census was the first national enumeration, it was closely linked to census-based conscription. A decisive shift occurred in 1861, following the establishment of a dedicated Department of Statistics, when Greece conducted its first scientifically organized population census, together with an Agricultural Census. The systematic compilation and publication of the results marked the emergence of population censuses as instruments of official statistical knowledge, rather than solely as tools serving immediate administrative purposes (Hellenic Statistical Authority, 2021; Kingdom of Greece, Ministry of the Interior, 1862).

The institutional reforms described above were accompanied by gradual improvements in the scope and organization of statistical production. The Directorate established in 1913 was explicitly responsible for collecting, verifying, processing and publishing national statistics. The centralization of statistical work in 1925 promoted methodological uniformity, staff specialization and the avoidance of duplicated statistical activities. In the following decades, the Statistical Service developed a regular program of censuses and annual and monthly publications covering population, agriculture, industry, trade, shipping, vital statistics and other fields (Hellenic Statistical Authority, 2021).

Further methodological milestones included the use of aerial photography for census mapping in the 1928 Population Census and the transition to a hybrid, digitally supported census model in 2021. The 2021 Population and Housing Census introduced electronic self-enumeration through a dedicated web application, supplemented by telephone or face-to-face interviews where necessary (Hellenic Statistical Authority, 2021; Benaki-Kyprioti and Kasapis, 2022).

Beyond enumerating the population, modern censuses also provide the sampling frame for a wide range of social, economic and environmental surveys and contribute to the development of statistical registers, thereby serving as core infrastructure for the broader system of official statistics (Hellenic Statistical Authority, 2021; Benaki-Kyprioti and Kasapis, 2022).

## 4. The Science of Probability and Statistics in Greece

Probability and statistics emerged as distinct academic disciplines in Greece between the late 1960s and 1982, when Law 1268/1982 abolished the former chair-based academic system (Hellenic Republic, 1982). Between 1968 and 1980, seven chairs in probability, statistics or closely related fields were established at the National and Kapodistrian University of Athens, the Higher School of Economics and Commercial Sciences, which was renamed the Athens University of Economics and Business in 1989 (Athens University of Economics and Business, 2021), Aristotle University of Thessaloniki, the University of Ioannina and the University of Patras. Six of the seven chairs were filled, while one of the two chairs at the University of Ioannina remained vacant until the abolition of the chair system in 1982. More precisely:

- **In 1968, Theophilos Cacoullos**, then an Associate Professor at New York University, assumed the Chair of Probability and Statistics at the National and Kapodistrian University of Athens. An Assistant Professorship position in Probability and Statistics, established a couple of years later, remained vacant until it was abolished in 1982. Cacoullos is well known for his research on variational inequalities, estimation, upper and lower bounds for the variance, and characterizations of distributions (Cacoullos, 1966, 1982; Cacoullos and Papageorgiou, 1983; Cacoullos and Papathanasiou, 1989; Cacoullos et al., 1994). He has 51 publications and was the editor of the Proceedings of an Advanced Study Institute on Discriminant Analysis and Applications.
- **In 1969, Demetrios Lambrakis**, who had received his PhD from the University of Aberdeen in 1966, was elected to the Chair of Probability and Statistics at the University of Ioannina. In 1973, he moved to the Higher School of Economics and Commercial Sciences, now the Athens University of Economics and Business, after being elected to a chair there. His research focused on experiments with mixtures and included an alternative to the simplex-lattice design (Lambrakis, 1969). His publication record includes five works in this field.

- **In 1972, George Roussas**, then a Professor at the University of Wisconsin–Madison, assumed the Chair of Applied Mathematics at the University of Patras. In 1984, he moved to the University of California, Davis, where he is now Distinguished Professor Emeritus of Statistics. His research made important contributions to nonparametric estimation and asymptotic distribution theory (Roussas, 1990; Roussas et al., 1992; Roussas, 1994; Tran et al., 1996), and he also authored influential books in mathematical statistics and probability (Roussas, 1997, 2008, 2024). His publication record includes 90 works. He supervised six PhD students and has 41 academic descendants. Five of his students are or have been university professors in Greece and/or the USA.

- **In 1976, Takis Papaioannou**, then an Assistant Professor at McGill University, was elected to the Chair of Probability and Statistics at the University of Ioannina. The second chair at the University of Ioannina remained vacant until the abolition of the chair system in 1982. In 1999, he moved to the University of Piraeus. His research made important contributions to parametric measures of information, ϕ-divergence measures, symmetry models for contingency tables, and entropy in actuarial science (Ferentinos and Papaioannou, 1981; Zografos et al., 1990; Kateri and Papaioannou, 1997; Sachlas and Papaioannou, 2014). His publication record includes 33 works. He supervised five PhD students and has six academic descendants. All five students are or have been university professors in Greece or Germany.

- **In 1976, Stratis Kounias**, then an Associate Professor at McGill University, was elected to the Third Chair of Mathematical Science at Aristotle University of Thessaloniki. In 1988, he moved to the National and Kapodistrian University of Athens. His research made important contributions to optimal designs, combinatorial probability, and Hadamard matrices (Kounias and Marin, 1976; Moyssiadis and Kounias, 1982; Kounias and Farmakis, 1988; Koukouvinos et al., 1991). His publication record includes 68 works. Based on the Mathematics Genealogy Project, Karagrigoriou (2009), and supplementary information provided by the authors, he supervised 12 PhD students and has more than 100 academic descendants. All 12 students are or have been university professors in Greece or the USA.
- **In 1980, Andreas Philippou**, then Professor at the Lebanese American University and Associate Professor at the American University of Beirut, was elected to the Second Chair of Applied Mathematics at the University of Patras. In 1988, he was invited to Cyprus to serve as Minister of Education from 1988 to 1990 and subsequently as a Member of Parliament from 1991 to 2001. He was granted leave of absence from the University of Patras during his public service, while continuing his research and the supervision of students. In 2001, he returned to the University of Patras. His pioneering research includes contributions to distributions of order k, Fibonacci polynomials of order k, runs and longest runs, and the reliability of consecutive systems (Philippou and Muwafi, 1982; Philippou et al., 1983; Philippou and Makri, 1985, 1986; Philippou, 1986, 1988; Philippou et al., 1988; Makri et al., 2007; Dafnis et al., 2010, 2019). His publication record includes 79 works. He supervised five PhD students and has 10 academic descendants. All five students are or have been university professors in Greece or Cyprus.

Of the six professors mentioned above, George Roussas was elected rector of his university, while Stratis Kounias and Andreas Philippou served as vice-rectors at their respective universities. The first five also wrote probability and statistics textbooks in Greek for their students.

Since 1982, more than 200 PhD holders in probability, statistics or related fields have been elected or promoted to academic positions across Greek public universities, reflecting the substantial expansion of these disciplines. Several of these professors have also served as rectors, vice-rectors, department chairs and deans. As of July 2026, twelve examples are mentioned below:

- **Professor Apostolos Burnetas** served as Chair of the Department of Mathematics at the National and Kapodistrian University of Athens.
- **Emeritus Professor George Donatos** served as Dean of the School of Law, Economics and Political Sciences at the National and Kapodistrian University of Athens. He currently serves as an External Member of the Board of Directors of the University of Piraeus.
- **Professor Alex Karagrigoriou** is the current Chair of the Department of Statistics and Insurance Science at the University of Piraeus. He also served as Chair of the Department of Statistics and Actuarial-Financial Mathematics at the University of the Aegean.
- **Professor Athanassios Katsis** is the current President of the Hellenic Statistical Authority (ELSTAT) and served as Rector of the University of the Peloponnese from 2017 to 2026.
- **Emeritus Professor Markos Koutras**, the guest of honor of the ΣtatMOD2025 Conference, served as Vice-Rector of the University of Piraeus.
- **Emeritus Professor Andreas Kyriakoussis** served as Rector of Harokopio University.
- **Professor John Panaretos** served as Vice-Rector of the University of Patras from 1988 to 1991 and as Deputy Minister of Education from 2009 to 2011.
- **Professor Stelios Psarakis** is the current Chair of the Department of Statistics at the Athens University of Economics and Business.
- **Professor Michael Sfakianakis** is the current Rector of the University of Piraeus.
- **Professor Vassilis GS Vasdekis** is the current Rector of the Athens University of Economics and Business. He also served as Vice-Rector from 2020 to 2024.
- **Professor Georgia Verropoulou** is currently a Vice-Rector of the University of Piraeus.
- **Professor Konstantinos Zografos** currently serves as an Internal Member of the Board of Directors of the University of Ioannina.

## 5. Probability and Statistics in Greece: Present and Future

The evolution of official statistics in Greece reflects a gradual transition from administrative population enumeration to an integrated statistical system based on specialized institutions, standardized methods and increasingly advanced data-collection technologies. In parallel, probability, statistics, biostatistics, actuarial science, demography and operations research have become established fields of teaching and research in Greek universities. The historical trajectory outlined in this article demonstrates both institutional continuity and increasing methodological specialization. The current generation is well positioned to contribute to future methodological and institutional developments in probability and statistics.

More recently, the modernization of official statistics has extended beyond data production to dissemination and public engagement using open formats, visual communication tools and initiatives aimed at improving statistical literacy and the accessibility of official data (Kasapis et al., 2026).